%% file: Article2024.tex
\documentclass{article}

\usepackage{arxiv}
\input{parameters}

\usepackage[utf8]{inputenc} 
\usepackage[T1]{fontenc}    
\usepackage{hyperref}       
\usepackage{url}            
\usepackage{enumitem}
\setlist[description]{leftmargin=0.25cm}
\usepackage{booktabs}       
\usepackage{amsfonts}       
\usepackage{nicefrac}       
\usepackage{algorithmic}
\usepackage{algorithm}
\usepackage{microtype}      
\usepackage{lipsum}
\usepackage{comment}
\usepackage{graphicx}
\usepackage{caption}
\usepackage{subcaption}
\usepackage{xcolor}
 
\graphicspath{ {./images/} }
\usepackage{amsmath}
\theoremstyle{definition}
\newtheorem{theorem}{Theorem}[section]

\newtheorem{corolaire}[theorem]{Corollary}

\newtheorem{definition}[theorem]{Definition}

\newtheorem{notation}[theorem]{Notations}
\newtheorem{problem}[theorem]{Problem}

\newtheorem{proposition}[theorem]{Proposition}
\newtheorem{remark}[theorem]{Remark}

\title{Sensitivity analysis of multiobjective linear programming from a geometric perspective}
\author{
 Mustapha Kaci\\
 Department of Mathematics, University of Oran Mohamed Boudiaf USTO-MB,\\ Oran, Algeria.\\
  {\color{blue}\texttt{kaci.mustapha.95@gmail.com} }
}

\begin{document}
\fontfamily{ptm}\selectfont
\maketitle

\begin{abstract}
Sensitivity analysis plays a crucial role in multiobjective linear programming (MOLP), where understanding the impact of parameter changes on efficient solutions is essential. This work builds upon and extends previous investigations. In this paper, we introduce a novel approach to sensitivity analysis in MOLP, designed to be computationally feasible for decision-makers studying the behavior of efficient solutions under perturbations of objective function coefficients in a two-dimensional variable space. This approach classifies all MOLP problems in \( S \subset \mathbb{R}^{2} \) by defining an equivalence relation that partitions the space of linear maps—comprising all sequences of linear forms on \(\mathbb{R}^2\) of length \( K \geq 2 \)—into a finite number of equivalence classes. Each equivalence class is associated with a unique subset of the boundary of \( S \). For any MOLP with \( K \) objective functions belonging to the same equivalence class, its set of efficient solutions corresponds to the associated subset of the boundary of \( S \). This approach is detailed and illustrated with a numerical example.
\end{abstract}


\section{Introduction}
\par

Sensitivity analysis in linear programming serves as a fundamental technique for assessing how variations in input data affect the optimal solution of a linear program. This analysis delves into the responsiveness of the optimal solution to changes in the coefficients of the objective function and the constraints. By determining the permissible intervals for these parameters, sensitivity analysis provides insights into the robustness and flexibility of the optimal solution. This capability is particularly crucial in dynamic environments where data can fluctuate due to market conditions, resource availability, or other external factors. Through sensitivity analysis, decision-makers can identify which parameters have the most significant impact on the solution, allowing for more resilient and adaptive planning. This approach not only enhances the reliability of the optimization model but also supports proactive management strategies by foreseeing potential disruptions and enabling timely adjustments \cite{dantzig1963linear}.

\par Multicriteria optimization is essential in fields such as engineering, economics, and management, where decisions must balance various competing interests by maximizing or minimizing multiple, often conflicting criteria simultaneously. Zeleny developed the concept of compromise programming and compiled a selected bibliography of works related to multicriteria decision-making \cite{zeleny1973compromise, zeleny1973bibliography}. To address these challenges, several approaches have been developed. One common method is weighting objective functions, where each objective is assigned a weight reflecting its relative importance \cite{int1}. This approach has been adopted in this study. Other techniques include the constraint method, which transforms some objectives into constraints \cite{int2}, and goal programming, which focuses on achieving specific target levels for each objective, minimizing deviations from these targets \cite{int3}. Additionally, the fuzzy programming approach incorporates uncertainty and imprecision in the objective functions and constraints to find a solution that aligns with the decision-maker's preferences in a more flexible manner \cite{int5, int4}. For a more comprehensive survey of methods, including these approaches, refer to \cite{KORNBLUTH1973193, roy1971problems, Spronk1981}. Several algorithms are available for solving all extreme non-dominated solutions: Yu and Zeleny \cite{yu1973set}, Yu \cite{yu1974cone}, Steuer \cite{steuer1974adex}, Philip \cite{philip1972algorithms}, Charnes and Cooper \cite{charnes1961management}, Evans and Steuer \cite{evans1973revised}.

\subsection{Sensitivity Analysis in Linear Programming}

Sensitivity analysis in linear programming, as introduced by Bradley, Hax, and Magnanti \cite{bradley1977applied}, involves examining how variations in the data impact the optimal solution. Traditional sensitivity analysis considers changes in the objective function coefficients and the right-hand side values, assessing the ranges within which these variations do not alter the optimal basis. This analysis helps determine shadow prices and reduced costs, which remain constant within specific ranges. By adding a parameter to the coefficients and using the final simplex tableau, one can determine the interval within which the coefficients can vary without changing the optimal solution \cite{bradley1977applied,dantzig1963linear,murty1976linear}. Kaci and Radjef \cite{kaci2022} introduced a  geometric approach to sensitivity analysis in linear programming. Their method utilizes concepts from affine geometry to present a fresh formulation of the sensitivity analysis problems. Specifically, they represent the coefficient vector of the objective function in polar coordinates, identifying the angles at which the solution remains invariant to changes.

For MOLP, sensitivity analysis becomes more complex due to the presence of multiple conflicting objective functions. A recently derived multicriteria simplex method is employed to explore fundamental properties within the decomposition of parametric space \cite{yu1975set}. This method introduces a novel type of parametric space that naturally emerges from its formulation. The study utilizes two computational approaches: an indirect algebraic method, which identifies the set of all nondominated extreme points, and a direct geometric decomposition method akin to the approach discussed by Gal and Nedoma \cite{gal1972multiparametric}. In addressing MOLP problems, Steuer developed a methodology that uses interval weights for each objective, rather than fixed weights \cite{steuer1976multiple}. This approach contrasts with traditional methods that typically yield a single efficient extreme point. When interval weights are applied, a cluster of efficient extreme points is generated. This cluster allows decision-makers to qualitatively identify the solution that offers the greatest utility, potentially simplifying the decision process by providing a range of efficient solutions close to the optimal point. The methodology converts the problem into an equivalent vector-maximum problem, utilizing algorithms designed for such problems to identify the subset of efficient extreme points corresponding to the specified interval weights. In their pioneering work, Richard E. Wendell introduces an innovative method for sensitivity analysis that permits independent and simultaneous variations in multiple parameters. This approach calculates the maximum tolerance percentages for both right-hand-side terms and objective function coefficients, ensuring the optimal basis remains unchanged within these tolerance ranges \cite{wendell1984tolerance,wen}.

\subsection{Advances and Contributions in Sensitivity Analysis}

In their paper \cite{kaci2024}, building on their previous work \cite{kaci2022}, Kaci and Radjef propose an approach to solve MOLP problems that avoids computing individual objective function optima. They partition the space of linear forms using a fixed convex polygonal subset of \(\mathbb{R}^2\) and an equivalence relation, ensuring elements within each equivalence class share the same optimal solution.

This study introduces a novel geometric approach to sensitivity analysis for MOLP problems with two decision variables. Our method visualizes the weighting of the objective function as a rotation of the graph of one of its components. This rotation occurs between the two extreme rays of the gradient cone, generating all efficient solutions and significantly simplifying calculations. This approach demonstrates the equivalence between MOLP and two objective linear programming (TOLP) problems, facilitating their classification.

\subsection{Paper Structure}
\par This paper is organized as follows: Section \ref{section 2} introduces the preliminary concepts essential for understanding the problem. Section \ref{section 3} details the formulation of the sensitivity analysis problem. Section \ref{section 4} explores the relationship between the two extreme rays and the extremal points of the efficient solution set of the MOLP problem,  and presents Algorithm \ref{ama} to construct the solution set for the sensitivity analysis problem. Section \ref{section 5} includes an illustrative numerical example. The discussion is presented in Section \ref{section 6}, and the paper concludes with Section \ref{section 7}.

\section{Preliminary}\label{section 2}
\subsection{Sensitivity analysis for a mono-objective inear programming problem}
In this section, we present the foundational tools necessary to establish the results of this paper. It is essential to note the close connection between this work and the study by Kaci and Radjef \cite{kaci2022}. Therefore, we recommend readers familiarize themselves with \cite{kaci2022} before delving into this paper.

Throughout our analysis, \( F^{0} \) will represent a function of two variables with values in \(\mathbb{R}^{K}\), where \( K \geq 2 \). Each \( f^{0}_{k} \) for \( k = 1, \ldots, K \) denotes a linear form defined on \(\mathbb{R}^{2}\) and taking values in \(\mathbb{R}\). Specifically, \( F^{0}(x) \) is defined as:
\[
F^{0}(x) = \left(f^{0}_{1}(x), f^{0}_{2}(x), \ldots, f^{0}_{K}(x)\right),\quad \text{for all} \; x\in\mathbb{R}^{2}.
\]
Consider the following initial mono-objective linear programming (MO-OLP) problem:
\begin{equation}\label{Pbm0}
\max_{x \in S} \hspace{0.15cm} f^{0}_{k}(x) = c_{0k}^{t}x,
\end{equation}
where 
\[ 
S := \left\{ x \in \mathbb{R}^2 \hspace{0.17cm}:\hspace{0.17cm} Ax \leq b, \, x \geq 0 \right\}, 
\]
with
\[
A=\left(\begin{array}{cc}
a_{11} & a_{12} \\
a_{21} & a_{22} \\
\vdots & \vdots \\
a_{m1} & a_{m2}
\end{array}\right), \quad
b=\left(\begin{array}{c}
b_{1} \\
b_{2} \\
\vdots \\
b_{m}
\end{array}\right), \quad
c_{0k}=\left(\begin{array}{c}
c_{1}^{0k} \\
c_{2}^{0k}
\end{array}\right), \quad
x=\left(\begin{array}{c}
x_{1} \\
x_{2}
\end{array}\right),
\]
where \( c_{0k}^t \) denotes the transpose of \( c_{0k} \), and \( b_i \), \( a_{i1} \), \( a_{i2} \), \( c_{1}^{0k} \), and \( c_{2}^{0k} \) for \( i = 1, \ldots, m \) and \( k = 1, \ldots, K \) are constants. Here, \( m \) represents the number of linear constraints.

\begin{notation}
Let \( S\subset\mathbb{R}^2  \) be a closed and bounded convex polygon.
\begin{itemize}
    \item The boundary of \( S \) is denoted by \( \partial S \).

    \item The vector space of linear forms on \( \mathbb{R}^{2} \), denoted \( \mathcal{L}(\mathbb{R}^{2}) \), consists of functions \( g: \mathbb{R}^{2} \rightarrow \mathbb{R} \) defined by \( g(x_{1}, x_{2}) = c_{1}x_{1} + c_{2}x_{2} \), where \( c_{1}, c_{2} \in \mathbb{R} \).

    \item 
The set of points \( y \in S \) such that \( f^{0}_{k}(y) \geq f^{0}_{k}(x) \) for all \( x \in S \) is denoted as:
$$
\mathcal{P}\left(f^{0}_{k}\right) = \arg \max_{x \in S} f^{0}_{k}(x).
$$
    \item 
The set of vertices is denoted by \( \mathcal{V}(S) \). Since \( S \) is a closed and bounded polygon in \(\mathbb{R}^{2}\), it has a finite number of vertices (extreme points). Let \( v^{1}, \ldots, v^{|\mathcal{V}(S)|} \) represent the vertices of \( S \), numbered counterclockwise, where \(|\mathcal{V}(S)|\) is the total number of vertices.
    \item Denote the following index sets: \(\mathcal{K} = \{1, \ldots, K\}\) and \(V(S) = \left\{1, \ldots, |\mathcal{V}(S)|\right\}\).
    \item The inner product in \(\mathbb{R}^2\) is denoted by \(\left\langle \cdot, \cdot \right\rangle\). Then, \( f_{k}^{0} \) for \( k \in \mathcal{K} \) will be denoted by 
\[ 
f_{k}^{0} = \left\langle c_{0k}^{t}, \cdot \right\rangle.
\]
\end{itemize}
\end{notation}

\begin{problem}[Sensitivity analysis for MO-OLP problem]\label{S2022}
Let \( k \in \mathcal{K} \). Find  all linear forms $g\in\mathcal{L}\left(\mathbb{R}^{2}\right)$ such that:
$$
\mathcal{P}\left(f^{0}_{k}\right)=\mathcal{P}\left(g\right).
$$
\end{problem}
\begin{remark}$\;$\\
\begin{enumerate}
   \item Our focus hereafter lies in conducting sensitivity analysis for a MOLP problem. Hence, we introduce the initial problem (\ref{Pbm2}), where the objective function is given by \( F^{0} \). This leads us to manipulate its components \( f_{k}^{0} \), which motivates their use in problems (\ref{Pbm0}) and \ref{S2022}.
    \item 
In \cite{kaci2022}, the solutions to problem \ref{S2022} are initially established assuming positive coefficients \( c_{1}^{0k} \) and \( c_{2}^{0k} \) for \( f_{k}^{0} \). Remark 4.3 in the same reference extends this result to the broader case where these coefficients are real numbers. Here, we present the solution set to problem \ref{S2022} in its general form.
\end{enumerate}

\end{remark}

\begin{figure}
    \centering    \includegraphics[width=16.5cm,height=8cm]{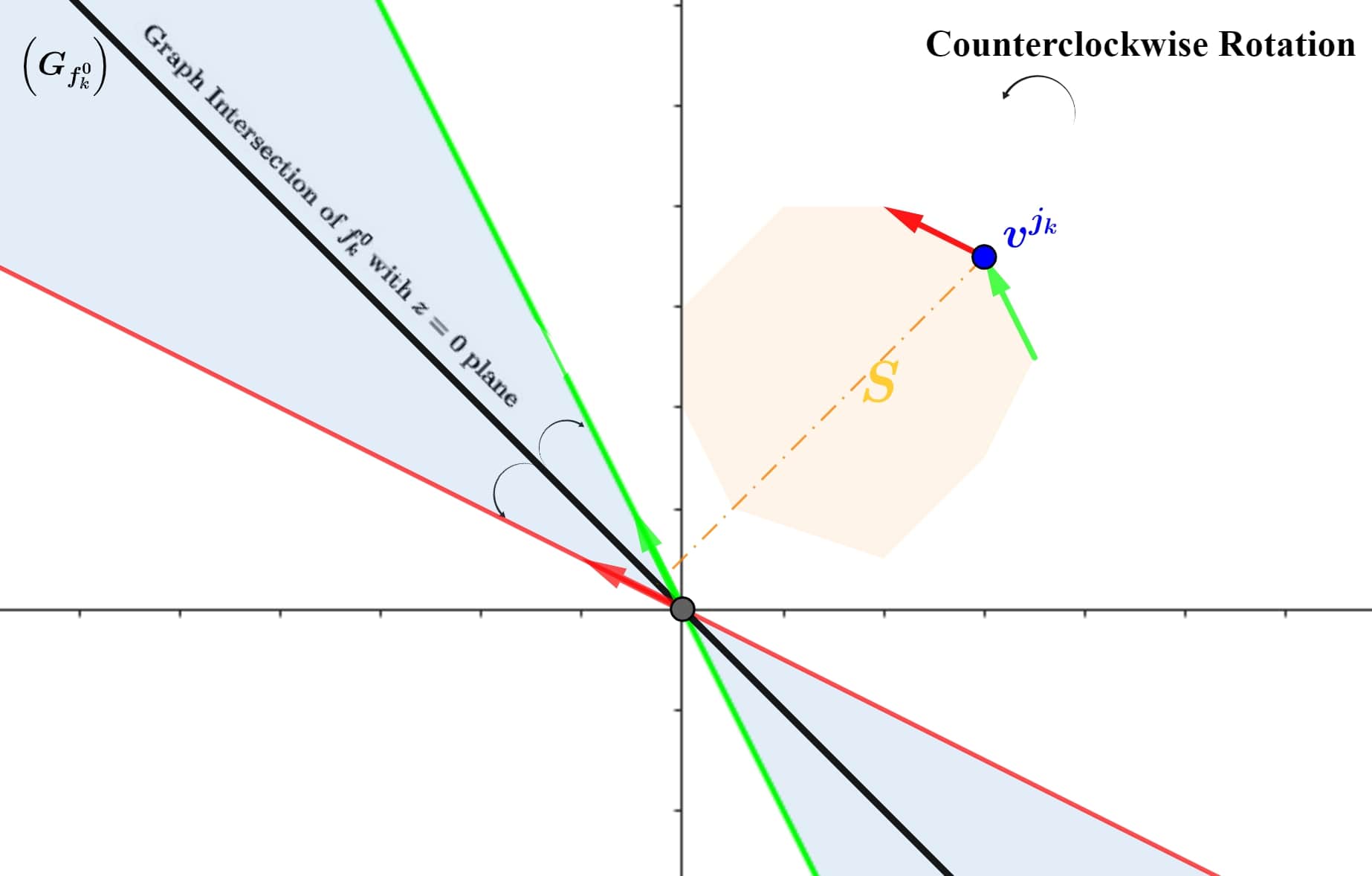}
        \caption{Geometric representation of the solution set \(\overline{f_{k}^{0}}\) given by Proposition \ref{SOL} for problem \ref{S2022}.\label{solp}}
\end{figure}

To tackle Problem \ref{S2022}, the authors adopted a geometric perspective, treating both the feasible region and the objective function as geometric objects. They demonstrated that the optimal solution \(v^{j_{k}}\) to Problem (\ref{Pbm0}) is located within \( \mathcal{V}(S) \) and maximizes the distance from its orthogonal projection onto the graph of the objective function. Furthermore, the solutions to Problem \ref{S2022} are all linear forms of \(\mathcal{L}\left(\mathbb{R}^{2}\right)\) whose intersection with the plane \(z=0\) (\(G_{f_{k}^{0}}\), the black vector line in Figure \ref{solp}) lies within the cone defined by the lines parallel to the vectors \(v^{j_{k}\mathcal{R}(j_{k}-1)}\) and \(v^{\mathcal{R}(j_{k}+1)j_{k}}\) — shown as green and red lines, respectively, in Figure \ref{solp}.

 This defined the range of angles within which the optimal solution is maintained. Details can be found in \cite{kaci2022}. 
 This rotational concept motivates the use of polar coordinates.

Let \( k \in \mathcal{K} \), and \( v^{j_{k}} \in \mathcal{V}(S) \), where \( j_{k} \in V(S) \) denotes the optimal solution of problem \eqref{Pbm0}. Let \( v^{j_{k}-1}, v^{\mathcal{M}(j_{k}+1)} \in V(S) \) such that \( v^{\mathcal{M}(j_{k}-1)}, v^{j_{k}}, v^{\mathcal{M}(j_{k}+1)} \) are successive corners of \( S \). Then, the following vectors are given in polar coordinates:
\begin{equation}\label{vb}
\begin{array}{lcccl}
   v^{j_{k}\mathcal{M}(j_{k}-1)} &=& v^{j_{k}} - v^{\mathcal{M}(j_{k}-1)} &= &  r_{1}(k)\left( \cos(\theta_{1}(k)), \sin(\theta_{1}(k)) \right) \vspace{0.2cm}\\
   v^{\mathcal{M}(j_{k}+1)j_{k}} &= &v^{\mathcal{M}(j_{k}+1)} - v^{j_{k}} &= & r_{2}(k) \left( \cos(\theta_{2}(k)), \sin(\theta_{2}(k)) \right)\vspace{0.2cm}\\
c_{0k}^{t} &= & \left(c_{1}^{0k},c_{2}^{0k}\right)&=& r(k) \left( \cos(\phi(k)), \sin(\phi(k)) \right),
\end{array} 
\end{equation}
where \( r(k), r_{1}(k), r_{2}(k) > 0 \), \( \phi(k), \theta_{1}(k), \theta_{2}(k) \in [0, 2\pi[ \), and for all integer \( \alpha \geq 0 \),  \( \mathcal{M}(\alpha) \) is defined as follows:
\[ 
\mathcal{M}(\alpha) = \begin{cases} 
\alpha \mod |\mathcal{V}(S)| & \text{if } \alpha \text{ is not a multiple of } |\mathcal{V}(S)|, \\
|\mathcal{V}(S)| & \text{otherwise}.
\end{cases} 
\]
\begin{remark}
For \( j \in V(S) \), \(\mathcal{M}(j) \) denotes the new index of \( v^{j} \) when \( j \) exceeds the size of \( \mathcal{V}(S) \).
\end{remark}
\begin{proposition}[\cite{kaci2022}]\label{SOL}
The solution set to problem \ref{S2022} is defined by:
\begin{equation}
 \overline{f_{k}^{0}}=\left\{g\in\mathcal{L}\left(\mathbb{R}^{2}\right)\hspace{0.1cm} : \hspace{0.1cm}g(x) = \left\langle(\cos(\omega),\sin(\omega)),x\right\rangle, \hspace{0.15cm} r > 0, \hspace{0.15cm} \omega \in [0, 2\pi[,\hspace{0.15cm}\theta_1(k) < \omega + \frac{\pi}{2} < \theta_2(k)\right\}.
\end{equation}
\end{proposition}

\subsection{Efficient solution set of a MOLP problem}

Consider the folowing Iinitial MOLP problem:
\begin{equation}\label{Pbm2}
      \underset{x\in S}\max\hspace{0.15cm} F^{0}(x)= \underset{x\in S}\max\hspace{0.15cm} \left(f^{0}_1(x),f^{0}_2(x),\cdots,f^{0}_K(x)\right),       
\end{equation}

\begin{definition}
    
A feasible solution \( v \in S \) is a non-dominated solution, Pareto optimal solution, or efficient solution for the problem (\ref{Pbm2}) if there does not exist another feasible solution \( x \in S \) such that:
$$
\left\{
\begin{array}{ll}
f^{0}_{k}(v) \leq f^{0}_{k}(x) & \text{for all } k \in \mathcal{K} \vspace{0.15cm}\\
f_{k_{0}}(v) < f_{k_{0}}(x) & \text{for some } k_{0} \in \mathcal{K}.
\end{array}
\right.
$$

\begin{notation}
     The set of efficient solutions of problem (\ref{Pbm2}) is denoted by \( \mathcal{S}(F^{0}) \).
\end{notation}   
\end{definition}

\begin{definition}[Adjacency and connectivity of vertices in a polygon\label{ad}]$\;$\\
\begin{enumerate}
    \item Two vertices of a polygon are said to be adjacent if they are connected by an edge of the polygon.
    \item Two vertices \( \widetilde{v} \) and \( \hat{v} \) in \( S \) are considered connected if there exists a sequence of vertices \( v^1 = \widetilde{v}, v^2, \ldots, v^{\hat{j}} = \hat{v} \), with \( \hat{j} \geq 2 \), such that each consecutive pair \( (v^i, v^{i+1}) \) for \( i = 1, 2, \ldots,  \hat{j} - 1 \) is adjacent.
    \item A subset \( \mathcal{V} \subseteq \mathcal{V}(S) \) is connected if every pair of elements in \( \mathcal{V} \) is connected.
\end{enumerate}
\end{definition}
\begin{theorem}[Theorem 3.7 \cite{yu1973set}\label{conn}]
     \( \mathcal{S}(F^{0})\) is connected.
\end{theorem}

\begin{definition}
    \( x \in S \) is an extremal point of \(  \mathcal{S}(F^{0}) \) if it belongs to \( \partial \mathcal{S}(F^{0})  \cap  \mathcal{V}(S)\). In other words, \( x \) is extremal point if there exists \( y \in \mathcal{V}(S) \cap (S \setminus \mathcal{S}(F^{0})) \) adjacent to \( x \). For an example, see Figure \ref{exp}.
\end{definition}

\section{Problem formulation}\label{section 3}

The problem under study involves performing a sensitivity analysis on the set \(\mathcal{S}(F^{0})\) of efficient solutions of the MOLP problem (\ref{Pbm2}). This analysis requires determining the set \(\mathcal{A}(F^{0}) \subseteq \mathcal{L}\) of all linear mappings \(G \in \mathcal{L}\) that satisfy:
\begin{equation}\label{san}
   \mathcal{S}(F^{0}) = \mathcal{S}(G), \quad \text{for all} \quad G \in \mathcal{A}(F^{0}),
\end{equation}
where 
\[
\mathcal{L} = \bigcup_{K \geq 2} \mathcal{L}\left(\mathbb{R}^{2}\right)^{K}.
\]
Consider the binary relation
\[
\forall\ G, H \in \mathcal{L}, \quad G \;\mathcal{R}_{S}\; H \iff \mathcal{S}(G) = \mathcal{S}(H).
\]
The set \(\mathcal{A}(F^{0})\) represents the equivalence class of \(F^{0}\) under the relation \(\mathcal{R}_{S}\). The set of all such equivalence classes is denoted by 
\[
\mathcal{A} := \mathcal{L}/\mathcal{R}_{S} = \bigcup_{F^{0} \in \mathcal{L}} \mathcal{A}(F^{0}).
\]
\begin{notation}
Let \( F^{0} \in \mathcal{L} \). We denote the set of efficient extreme points contained in \(\mathcal{S}(F^{0})\) by \(\mathcal{VS}(F^{0})\), defined as \(\mathcal{VS}(F^{0}) = \mathcal{S}(F^{0}) \cap \mathcal{V}(S)\).

For \( F^{0} \in \mathcal{L}(\mathbb{R}^{2})^{K} \) and given \( j_{k_{1}}, j_{k_{2}} \in V(S) \) such that \( j_{k_{2}} \geq j_{k_{1}} \) for some \( k_{1}, k_{2} \in \mathcal{K} \), the set of efficient extreme points is expressed as:
\begin{equation}\label{vs}
    \mathcal{VS}(F^{0}) = \left\{ v^{j_{k_{1}}}, v^{j_{k_{1}}+1}, \ldots, v^{j_{k_{2}}} \right\}.
\end{equation}
\end{notation}

\section{Building the set 
$\mathcal{A}\left(F^{0}\right)$}\label{section 4}
\begin{figure}
    \centering    \includegraphics[width=16cm,height=10cm]{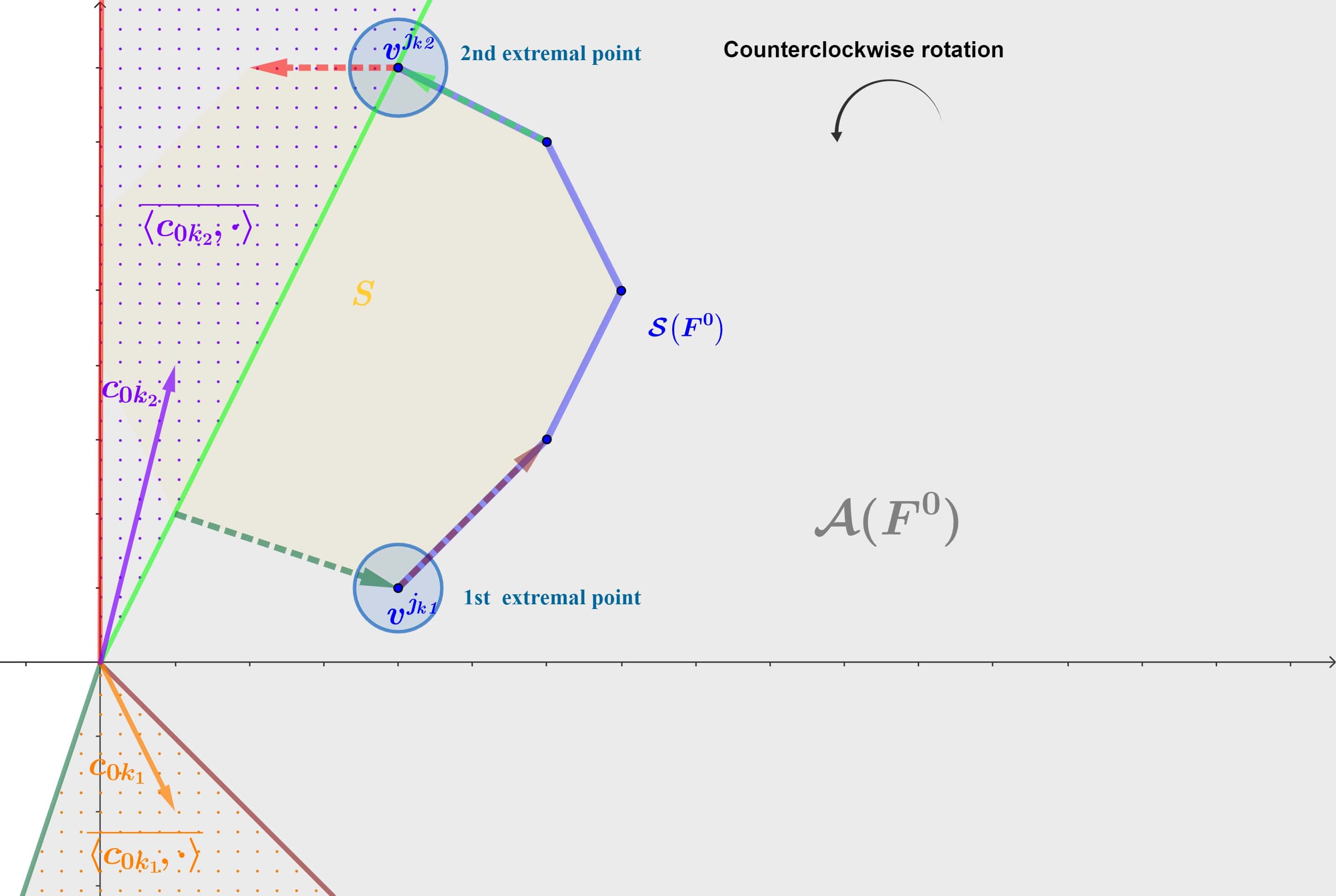}
        \caption{Geometric representation of the solution set \(\mathcal{A}(F^{0})\)  of the problem described in Section \ref{section 3}. This figure illustrates the example in Section \ref{section 5} with $k_{1}=5$ and $k_{2}=6$.\label{GH}}
\end{figure}

In this section, we detail the construction of the set $\mathcal{A}\left(F^{0}\right)$. Figure \ref{GH} provides a graphical representation of this process, which will aid in visualizing the concepts discussed. We encourage readers to refer to the figure for a clearer understanding of the construction steps.

The construction of \(\mathcal{A}\left(F^{0}\right)\) is driven by key considerations in two dimensions (\(S \subset \mathbb{R}^{2}\)). The assumption that \(S\) is a closed, bounded set in \(\mathbb{R}^{2}\) implies that \(\partial S\) forms a polygon. Additionally, \(\mathcal{S}(F^0)\) represents a continuous polygonal curve with only two extremal points.
\(\mathcal{S}(F^0)\) is given by:
$$
\mathcal{S}(F^0)=\bigcup_{d\in\mathcal{G}(F^0)}{\left(\mathcal{S}\left\langle d^{t}, \cdot \right\rangle\right)},
$$
where $\mathcal{G}(F^0)$ is the gradient cone defined by:
    $$
        \mathcal{G}(F^0)=\left\{wd \in \mathbb{R}^{2}, w \geq 0 \; : \; d = \sum_{k=1}^{K}{\lambda_{k}c_{0k}}, \; \lambda\in\Lambda(K)^{0}\right\},
    $$
with
$$
\Lambda(K)^{0}=\left\{\lambda\in\mathbb{R}^{K} \; : \;  \lambda_{k} \in [0,1],  \; \; \sum_{k=1}^{K} \lambda_{k} = 1\right\}.
$$
\begin{remark}\label{lf} Since \(\mathcal{G}(F^0) \subset \mathbb{R}^{2}\), \(\mathcal{G}(F^0)\) is a cone generated by only two extreme rays, i.e., there exist \(k_{1}, k_{2} \in \mathcal{K}\) such that
        \[
        \mathcal{G}(F^0) = \left\{wd \in \mathbb{R}^{2}, w \geq 0 \; : \; d = \delta c_{0k_{1}} + (1 - \delta) c_{0k_{2}}, \; \delta \in [0,1] \right\},
        \]
        where \(c_{0k_{1}}\) and \(c_{0k_{2}}\) generate the two extreme rays of \(\mathcal{G}(F^0)\); see Figure \ref{gco}.
\end{remark}

\subsection{Relationship between two extreme rays of \(\mathcal{G}(F^{0})\) and the extremal points of $\mathcal{S}(F^{0})$}
Consider the weighted objective function problem for \(F^{0}\), described as follows:
\begin{equation}\label{PND}
    \underset{x \in S}{\max} \, F^{0}_{\lambda}(x) = \underset{x \in S}{\max} \, 
 \lambda^{t} F^{0}(x)  \; \; \;
 \text{with} \; \; \; \lambda\in\Lambda(K)^{0},
\end{equation}
Since for all $\lambda \in [0,1]^{K}$, \(F^{0}_{\lambda}\) attains its maximum value at an extreme point when the distance between this point and its orthogonal projection onto the line vector defined by the intersection of the graph of \(F^{0}_{\lambda}\) with the plane \(z=0\) is maximized (see Proposition 4.4 in \cite{kaci2022}), we will focus exclusively on rotations of the plane \(z=0\). This implies that finding \(\mathcal{A}(F^{0})\) amounts to working in polar coordinates and determining the interval within which the rotation angle is allowed to vary so that \(\mathcal{S}(F^{0})\) remains unchanged.
\begin{remark}\label{sr}
    Throughout this work, all angles are defined within the interval \([0, 2\pi[\). Consequently, all rotations are considered to be in the positive direction (counterclockwise).
\end{remark}

\begin{notation}
Let us express the two generators of the extremal rays in polar coordinates:
$$
c_{0k_{i}}^{t} = r(k_{i}) \left( \cos(\phi(k_{i})), \sin(\phi(k_{i})) \right), \quad \text{for} \quad i = 1, 2,
$$
where $r(k_{1}), r(k_{2}) > 0$, and $\phi(k_{1}), \phi(k_{2}) \in [0,2\pi[$, and $\phi(k_{2}) \geq \phi(k_{1})$.
\end{notation}
\begin{theorem}\label{lem}
The map \(\delta \mapsto d_{\delta} = \delta c_{0k_{1}} + (1 - \delta) c_{0k_{2}}\) for all \(\delta \in [0,1]\). \( d_{\delta} \) is a rotation composed with a homothety, i.e., there exists a rotation \(R_{\theta}\) of angle \(\theta\in[0,2\pi[\), and a homothety \(H_{\alpha}\) of ratio \(\alpha\in\mathbb{R}^{+}\), such that 
$$
d_{\delta} = H_{\alpha}(R_{\theta}(c_{0k_{1}})).
$$
\end{theorem}
\begin{proof}
Let \(\delta \in [0,1]\). The proof consists of determining, for all \(\delta\), the rotation and the homothety as functions of \(\delta\), \(c_{0k_{1}}\), and \(c_{0k_{2}}\), that is, finding \(\alpha \in \mathbb{R}\) and \(\theta \in [0,2\pi[\) such that
$$
d_{\delta} = H_{\alpha}(R_{\theta}(c_{0k_{1}})) = \delta c_{0k_{1}} + (1 - \delta) c_{0k_{2}}.
$$
That is to say,
\begin{equation}\label{th}
\alpha r\left(k_{1}\right)\left(\begin{array}{cr}
\cos(\theta) & -\sin(\theta) \\
\sin(\theta) & \cos(\theta)
\end{array}\right) \left(\begin{array}{c}
\cos(\phi(k_{1})) \\
\sin(\phi(k_{1}))
\end{array}\right) = \delta \left(\begin{array}{c}
c_{1}^{0k_{1}} \\
c_{2}^{0k_{1}}
\end{array}\right) + (1 - \delta) \left(\begin{array}{c}
c_{1}^{0k_{2}} \\
c_{2}^{0k_{2}}
\end{array}\right),
\end{equation}
with
$$
R_{\theta} = \left(\begin{array}{cr}
\cos(\theta) & -\sin(\theta) \\
\sin(\theta) & \cos(\theta)
\end{array}\right), \quad
c_{0k_{1}} = \left(\begin{array}{c}
c_{1}^{0k_{1}} \\
c_{2}^{0k_{1}}
\end{array}\right) = r\left(k_{1}\right) \left(\begin{array}{c}
\cos(\phi(k_{1})) \\
\sin(\phi(k_{1}))
\end{array}\right),
$$
where \(r\left(k_{1}\right) > 0\), \(\phi(k_{1}) \in [0,2\pi[\). Then, formula (\ref{th}) becomes:
$$
\left\{\begin{aligned}
\alpha r\left(k_{1}\right)\left[\cos(\theta)\cos(\phi(k_{1})) - \sin(\theta)\sin(\phi(k_{1}))\right] &= \alpha r\left(k_{1}\right)\cos\left(\phi(k_{1}) + \theta\right) = \delta c_{1}^{0k_{1}} + (1 - \delta) c_{1}^{0k_{2}}, \vspace{0.25cm}\\
\alpha r\left(k_{1}\right)\left[\sin(\theta)\cos(\phi(k_{1})) + \cos(\theta)\sin(\phi(k_{1}))\right] &= \alpha r\left(k_{1}\right)\sin\left(\phi(k_{1}) + \theta\right) = \delta c_{2}^{0k_{1}} + (1 - \delta) c_{2}^{0k_{2}},
\end{aligned}\right.
$$
Therefore, we have on one hand
\begin{equation}\label{exx}
\left\{\begin{array}{l}
 \left[ \alpha r\left(k_{1}\right) \cos(\phi(k_{1}) + \theta) \right]^2 = \left( \delta c_{1}^{0k_{1}} + (1 - \delta) c_{1}^{0k_{2}} \right)^2, \vspace{0.25cm}\\
 \left[\alpha r\left(k_{1}\right) \sin(\phi(k_{1}) - \theta) \right]^2 = \left( \delta c_{2}^{0k_{1}} + (1 - \delta) c_{2}^{0k_{2}} \right)^2.
\end{array}\right.
\end{equation}
By adding these two expressions (\ref{exx}), we obtain:
\[
\left[\alpha r\left(k_{1}\right) \cos(\phi(k_{1}) + \theta) \right]^2 + \left[\alpha r\left(k_{1}\right) \sin(\phi(k_{1}) + \theta) \right]^2 = \left( \delta c_{1}^{0k_{1}} + (1 - \delta) c_{1}^{0k_{2}} \right)^2 + \left( \delta c_{2}^{0k_{1}} + (1 - \delta) c_{2}^{0k_{2}} \right)^2.
\]
Therefore,
\[ 
\alpha^2 r^2(\phi(k_{1})) \underbrace{\left[ \cos^2(\phi(k_{1}) + \theta) + \sin^2(\phi(k_{1}) + \theta) \right]}_{=1} = \left( \delta c_{1}^{0k_{1}} + (1 - \delta) c_{1}^{0k_{2}} \right)^2 + \left( \delta c_{2}^{0k_{1}} + (1 - \delta) c_{2}^{0k_{2}} \right)^2.
\]
This yields the expression for the ratio of the homothety \(\alpha\) as follows:
\begin{equation}\label{rat}
    \alpha = \frac{\sqrt{\left( \delta c_{1}^{0k_{1}} + (1 - \delta) c_{1}^{0k_{2}} \right)^2 + \left( \delta c_{2}^{0k_{1}} + (1 - \delta) c_{2}^{0k_{2}} \right)^2}}{r\left(k_{1}\right)}.
\end{equation}
On the other hand,
\begin{equation}\label{rra}
 \alpha r\left(k_{1}\right)\cos(\theta + \phi(k_{1}))  =  \delta c_{1}^{0k_{1}} + (1 - \delta) c_{1}^{0k_{2}}  \Longleftrightarrow  \theta = \arccos\left(\frac{ \delta c_{1}^{0k_{1}} + (1 - \delta) c_{1}^{0k_{2}} }{ \alpha r\left(k_{1}\right)}\right) - \phi(k_{1}).
\end{equation}
Finally, by substituting the value of 
 $\alpha$  obtained from formula (\ref{rat}) into formula (\ref{rra}), we obtain:
\begin{equation}\label{tet}
\theta = \arccos\left(\frac{ \delta c_{1}^{0k_{1}} + (1 - \delta) c_{1}^{0k_{2}} }{ \sqrt{\left( \delta c_{1}^{0k_{1}} + (1 - \delta) c_{1}^{0k_{2}} \right)^2 + \left( \delta c_{2}^{0k_{1}} + (1 - \delta) c_{2}^{0k_{2}} \right)^2}}\right) - \phi(k_{1}).    
\end{equation}
\end{proof}
  \begin{remark}\label{rem}     
Let \( d \in \mathcal{G}(F^0) \). Then, \( \mathcal{P}(\langle d^t, \cdot \rangle) = \mathcal{P}(\langle wd^t, \cdot \rangle) \) for any \( w > 0 \). Thus, to generate \( \mathcal{S}(F^0) \), we can use gradient vectors only with a constant norm equal to $\|c_{0k_{1}}\|$, where \(\|\cdot\|\) is the Euclidean norm in \(\mathbb{R}^2\). Specifically, \( \mathcal{S}(F^0) \) can be generated as follows:
\begin{equation}\label{fF}
\mathcal{S}(F^0) = \bigcup_{d \in \hat{\mathcal{G}}}{\mathcal{P}\left( \langle d^t, \cdot \rangle \right)}
\end{equation}
with
$$
\hat{\mathcal{G}}(F^0) = \left\{ d \in \mathcal{G}(F^0) \; : \; \|d\| = \|c_{0k_{1}}\| \right\}.
$$
 \end{remark} 
Remark \ref{rem} is extremely important for the subsequent discussion. Firstly, it clarifies in formula \ref{fF} that the norm of the gradient of the weighted function does not play a role in obtaining efficient solutions. Secondly, it allows us to exclude the homotheties experienced by \(F^{0}_{\lambda}\) when varying \(\lambda\), and to consider only rotations.

Assume, without loss of generality, that \(\|c_{0k_{1}}\| = \|c_{0k_{2}}\|\) and \(\alpha = 1\). This means that for all \(\delta \in [0, 1]\), \( d_{\delta} = R_{\theta}(c_{0k_{1}}) \), where \(\theta\) is related to \(\delta\) by formula (\ref{tet}). Consequently, the findings we establish for this specific scenario will be valid for the general case as well.
  
Theorem \ref{lem} and the remark \ref{rem} imply that, graphically, to generate \(\mathcal{S}(F^0)\), we start with a vector \(c_{0k_{1}}\), which lies on the first extremal ray. We then apply rotations of angle \(\theta \in I(k_{1},k_{2}):=\left[0, (\phi(k_{2}) - \phi(k_{1}))\right]\) to this vector. For each \(\theta\), we optimize the linear form \(\langle R_{\theta}(c_{0k_{1}})^t, \cdot \rangle\) to obtain all the elements of \(\mathcal{S}(F^0)\). This demonstrates the following corollary:
\begin{corolaire}\label{ree}
    \(\mathcal{S}(F^0)\) is given by
    \[
    \mathcal{S}(F^0) = \bigcup_{\theta \in I(k_{1}, k_{2})} \mathcal{P}\left( \langle R_{\theta}(c_{0k_{1}})^t, \cdot \rangle \right) = \bigcup_{\theta \in \hat{I}(k_{1}, k_{2})} \mathcal{P}\left( \langle R_{\theta}(c_{0k_{1}})^t, \cdot \rangle \right),
    \]
    where \(\hat{I}(k_{1}, k_{2})\) is defined by:
    \[
    \hat{I}(k_{1}, k_{2}) = \left] \theta_{1}(k_{1}) - \frac{\pi}{2}, \theta_{2}(k_{2}) - \frac{\pi}{2} \right[,
    \]
    and \(\theta_{1}(k_{1})\) and \(\theta_{2}(k_{2})\) are as defined in formula (\ref{pok}).
\end{corolaire}

An immediate question that arises from the above is whether the points obtained from gradient vectors lying on the two extremal rays are extremal points. The answer is provided by the following result.

\begin{corolaire}\label{re}
   Consider the two generators \(c_{0k_{1}}, c_{0k_{2}}\) of the two extreme rays of \(\hat{\mathcal{G}}(F^0)\) and the points \(v^{j_{k_{1}}}, v^{j_{k_{2}}} \in \mathcal{V}(S)\), such that:
    \[
   v^{j_{k_{i}}}\in\mathcal{P}\left( \langle c_{0k_{i}}^t, \cdot \rangle \right),\quad \text{for} \quad i=1,2.
    \]
    where \(j_{k_{1}}, j_{k_{2}} \in V(S)\)  and \(j_{k_{1}} \leq j_{k_{2}}\). Then, \( v^{j_{k_{1}}} \) and \( v^{j_{k_{2}}} \) are extremal points.
\end{corolaire}

\begin{proof}
For \( i = 1, 2 \) let 
\begin{equation}\label{pok}
\begin{array}{lcccl}
   v^{j_{k_{i}}\mathcal{M}(j_{k_{i}}-1)} &=& v^{j_{k_{i}}} - v^{\mathcal{M}(j_{k_{i}}-1)} &= &  r_{1}(k_{i})\left( \cos(\theta_{1}(k_{i})), \sin(\theta_{1}(k_{i})) \right) \vspace{0.2cm}\\
   v^{\mathcal{M}(j_{k_{i}}+1)j_{k_{i}}} &= &v^{\mathcal{M}(j_{k_{i}}+1)} - v^{j_{k_{i}}} &= & r_{2}(k_{i}) \left( \cos(\theta_{2}(k_{i})), \sin(\theta_{2}(k_{i})) \right)\vspace{0.2cm}\\
c_{0k_{i}}^{t} &= & \left(c_{1}^{0k_{i}},c_{2}^{0k_{i}}\right)&=& r(k_{i}) \left( \cos(\phi(k_{i})), \sin(\phi(k_{i})) \right),
\end{array} 
\end{equation}
where \( r(k_{i}), r_{1}(k_{i}), r_{2}(k_{i}) > 0 \) and \( \phi(k_{i}), \theta_{1}(k_{i}), \theta_{2}(k_{i}) \in [0, 2\pi[ \). Then, according to Corollary \ref{ree}, if we take \(\theta = \phi(k_{2}) - \phi(k_{1})\), we obtain
\[
   R_{\theta}(v^{j_{k_{1}}}) = v^{j_{k_{2}}}\in\mathcal{P}\left( \langle c_{0k_{2}}^t, \cdot \rangle\right) = \mathcal{P}\left( \langle R_{\theta}(c_{0k_{1}}^t), \cdot \rangle  \right).
\]
\end{proof}

\subsection{Characterization of the Set \(\mathcal{A}(F^{0})\)}
In this section, we provide a classification of MOLP problems in \(\mathbb{R}^{2}\) based on their solution sets. This classification is more general than the one presented in \cite{kaci2024}. Subsequently, we will deduce the set \(\mathcal{A}(F^{0})\).

We begin by exploring whether \(\mathcal{S}(F^0)\) can be obtained by solving a MOLP problem with only two objectives.
Specifically, we investigate whether every MOLP problem can be associated with a TOLP problem such that both problems yield the same set of solutions. The theorem presented next provides insight into this question.

\begin{theorem}[Equivalence of TOLP and MOLP problems]\label{fff}
For any \( F^{0} \in \mathcal{L} \), there exists \( H \in \mathcal{L}\left(\mathbb{R}^{2}\right)^{2} \) such that
\begin{equation}\label{Spb_20}
\mathcal{S}(H) = \mathcal{S}(F^{0}).
\end{equation}
Moreover, \( H \) is given by
\begin{equation}\label{H}
H = (h_{1}, h_{2}) = \left(\langle c_{0k_{1}}, \cdot \rangle, \langle c_{0k_{2}}, \cdot \rangle\right).
\end{equation}
Conversely, for any \( H \in \mathcal{L}\left(\mathbb{R}^{2}\right)^{2} \) as defined in formula (\ref{H}), there exists a corresponding \( G \in \mathcal{A}(F^{0}) \), where the order of components is irrelevant, such that equation (\ref{Spb_20}) holds. Additionally,  \(G\) can be expressed as:
\begin{equation}\label{ff}
G = \left( \langle c_{0k_{1}}^{t}, \cdot \rangle, \langle c_{0k_{2}}^{t}, \cdot \rangle, \langle R_{\theta_{1}}(c_{0k_{1}})^t, \cdot \rangle, \ldots, \langle R_{\theta_{K-2}}(c_{0k_{1}})^t, \cdot \rangle \right), \quad \text{with} \quad \theta_{1}, \ldots, \theta_{K-2} \in \hat{I}(k_{1},k_{2}).
\end{equation}
\end{theorem}

\begin{proof}
Let \( F^{0} \in \mathcal{L}\left(\mathbb{R}^{2}\right)^{K} \) for some \( K \geq 2 \). By Remarks \ref{lf} and \ref{sr}, the vectors \( c_{0k_{1}} \) and \( c_{0k_{2}} \) generate the gradient cone \(\mathcal{G}(F^{0})\), which implies that \(\mathcal{G}(H) = \mathcal{G}(F^{0})\). Therefore,
\[
\mathcal{S}(H) =\mathcal{S}(F^{0}).
\]
Conversely, since \(\theta_{1}, \ldots, \theta_{K-2} \in [0, \phi(k_{2}) - \phi(k_{1})]\), it follows that \(R_{\theta_{k}}(c_{0k_{1}})^t \in \mathcal{G}(F^{0})\) for $k=1,\ldots,K-2$. Furthermore, the two generators of \(\mathcal{G}(F^{0})\) are the gradients of \(\langle c_{0k_{1}}, \cdot \rangle\) and \(\langle c_{0k_{2}}, \cdot \rangle\), which are components of \(G\). Consequently, \(\mathcal{G}(G) = \mathcal{G}(F^{0})\), implying that \(G \in \mathcal{A}(F^{0})\). Additionally, since \(\mathcal{G}(G) = \mathcal{G}(H)\), it follows that \(\mathcal{S}(G) = \mathcal{S}(H)\).

\end{proof}

\begin{corolaire}\label{ffd}
$\mathcal{A}(F^{0})$ is given by:
  \begin{equation}\label{sets}
\mathcal{A}(F^{0}) = \bigcup_{K \geq 2} \left( \bigcup_{(\theta_{1}, \ldots, \theta_{K-2}) \in \hat{I}(k_{1},k_{2})^{K-2}} \left\{ \left( \langle c_{0k_{1}}^{t}, \cdot \rangle, \langle c_{0k_{2}}^{t}, \cdot \rangle, \langle R_{\theta_{1}}(c_{0k_{1}})^t, \cdot \rangle, \ldots, \langle R_{\theta_{K-2}}(c_{0k_{1}})^t, \cdot \rangle \right) \right\} \right).
\end{equation}  
\end{corolaire}
\begin{proof}
Theorem \ref{fff} ensures the existence of \( G \in \mathcal{A}(F^{0}) \) as defined in formula (\ref{ff}). Hence, it suffices to consider all elements of this form by taking all possible angles and any number of objectives \( K \geq 2 \).
\end{proof}

Now, let us review the key information about efficient solutions: their locations and the geometric form of \(\mathcal{S}(F^{0})\). The efficient solutions of MOLP problem form a continuous polygonal curve along \(\partial S\), where the vertices of this curve are the extreme efficient solutions. To determine all possible sets of efficient solutions on \(\partial S\) for any given MOLP problem, our next goal is to identify the subsets of efficient solutions that can be obtained from solving the problem (\ref{Pbm2}). 
\begin{notation}
Let \(\mathcal{NS}\) denote the set of all possible subsets of \(\partial S\) that can appear as efficient solution sets for an arbitrary MOLP problem. Therefore, the efficient solution set of any given MOLP problem belongs to \(\mathcal{NS}\).
\end{notation}

\begin{definition} 
The construction of \( \mathcal{NS} \) is given in the following steps:
\begin{enumerate}
    \item Let \(\mathcal{S}_{j^{0}}^{j}\) denote the polygonal curve formed by the \(j\) successive vertices of \(S\), starting from the \(j^{0}\)th vertex and continuing to the \((j^{0} + j - 1)\)th vertex. This is defined as:
    \begin{equation}\label{cur}
       \mathcal{S}_{j^{0}}^{j} = \bigcup_{l=1}^{j-1} \mathcal{H}\left\{v^{\mathcal{M}(j^{0}+l-1)}, v^{\mathcal{M}(j^{0}+l)}\right\},  
    \end{equation}    
    where for \( l = 1, \ldots, j-1 \),  \(\mathcal{H}\left\{v^{\mathcal{M}(j^{0}+l-1)}, v^{\mathcal{M}(j^{0}+l)}\right\}\) denotes the convex hull of the extreme points \( v^{\mathcal{M}(j^{0}+l-1)}\) and \(v^{\mathcal{M}(j^{0}+l)} \).
    \item Define the singletons containing the \(j^{0}\)th vertices as follows:
\[
\mathcal{S}^{1}_{j^{0}} = \left\{v^{j^{0}}\right\}.
\]
    \item Finally, \(\mathcal{NS}\) is defined as follows:
    \[
    \mathcal{NS} = \bigcup_{j^{0}=1}^{|V(S)|}{\left(\bigcup_{j=1}^{|V(S)|}{\mathcal{S}_{j^{0}}^{j}}\right)}.
    \]
\end{enumerate}
\end{definition}

\begin{remark}[Classification of MOLP problems]\label{CLA}
    From Formula (\ref{vs}), we see that \(\mathcal{S}(F^{0})\) can also be represented as \(\mathcal{S}_{j_{k_{1}}}^{j_{k_{2}}-j_{k_{1}}+1}\), where \(j_{k_{1}}, j_{k_{2}} \in V(S)\) and \(j_{k_{2}} \geq j_{k_{1}}\). 
    More generally, for every \(\mathcal{S}^{j}_{j^{0}} \in \mathcal{NS}\) where \(j^{0}\) and \(j\) are elements of \(V(S)\), there exists a unique class \(\mathcal{A}(F)\) with \(F \in \mathcal{L}\) such that:
\[
\mathcal{S}(G) = \mathcal{S}^{j}_{j^{0}}, \quad \text{if and only if} \quad G \in \mathcal{A}(F).
\]
\end{remark}

\begin{algorithm}[H]
\caption{Generate the solution set \(\mathcal{A}(F^{0})\) for the sensitivity analysis problem (\ref{san}).\label{ama}}

\bigskip

\textbf{Input:} \(K \geq 2\), \(\mathcal{K} = \{1, \ldots, K\}\), objective function \(F^{0} \in \mathcal{L}(\mathbb{R}^{2})^{K}\), feasible region \(S\), set of vertices \(\mathcal{V}(S)\), index set of vertices \(V(S)\), and consider the \(K\)-objectives problem (\ref{Pbm2}).\\
\textbf{Output:} \(\mathcal{A}(F^{0})\) solution to the sensitivity analysis problem (\ref{san}).

\bigskip

\begin{description}
\item[Part 1: Generate the efficient solution set \( \mathcal{S}\left(F^{0}\right) \)  of problem (\ref{Pbm2}).] $\;$\\
\begin{description}
    \item[Step 1:] Express the gradient vectors \(c_{0k}\) for \(k \in \mathcal{K}\) in polar coordinates, as follows:
\[
c_{0k} = r(k) \left( \cos(\phi(k)), \sin(\phi(k)) \right), \quad \text{for} \quad  k \in \mathcal{K}, \; \;  \text{and} \; \; \phi(k) \in [0, 2\pi[.
\]
    \item[Step 2:] Define the two generators \(c_{0k_{1}}\), \(c_{0k_{2}}\) for \(k_{1}, k_{2} \in \mathcal{K}\) of the two extreme rays of \(\mathcal{G}(F^{0})\), such that 
\[
\phi(k_{1}) = \min_{k \in \mathcal{K}} \phi(k), \quad \text{and} \quad \phi(k_{2}) = \max_{k \in \mathcal{K}} \phi(k).
\]
    \item[Step 3:] Apply Corollary \ref{re}, and get the two extremal efficient solutions \(v^{j_{k_{1}}}\) and \(v^{j_{k_{2}}}\), see Figure \ref{exp}. Then,
\[
\mathcal{VS}(F^{0}) = \{v^{j} \; : \; j = j_{k_{1}}, \ldots, j_{k_{2}}\}.
\]
    \item[Step 4:] Use formula (\ref{cur}) to obtain \(\mathcal{S}(F^{0})\) by taking \(j^{0} = j_{k_{1}}\) and \(j = j_{k_{2}} - j_{k_{1}} + 1\), as follows:
\[
\mathcal{S}(F^{0}) := \mathcal{S}_{j_{k_{1}}}^{j_{k_{2}} - j_{k_{1}} + 1} = \bigcup_{l = 1}^{j_{k_{2}} - j_{k_{1}}} \mathcal{H}\{v^{\mathcal{R}(j_{k_{1}} + l - 1)}, v^{\mathcal{R}(j_{k_{1}} + l)}\}.
\]
\end{description} 
\item[Part 2: Generate the solution set \(\mathcal{A}(F^{0})\) for the sensitivity analysis problem (\ref{san}).] $\;$\\
\begin{description}
    \item[Step 5:] Calculate \(\theta_{1}(k_{1})\) and \(\theta_{2}(k_{2})\) using formulas (\ref{pok}).
    \item[Step 6:] Obtain the solution set  \(\mathcal{A}(F^{0})\) to the sensitivity analysis problem (\ref{san}) using Corollary \ref{ffd}.
\end{description} 
\end{description} 
\end{algorithm}

\section{Numerical example}\label{section 5}

Consider the following six-objective linear programming problem:
\begin{equation}\label{exemple2}
\underset{x \in S}{\max} \, F^{0}(x)=\left(f_{1}^{0}(x),\hspace{0.17cm}f_{2}^{0}(x),\hspace{0.17cm}f_{3}^{0}(x),\hspace{0.17cm}f_{4}^{0}(x),\hspace{0.17cm}f_{5}^{0}(x),\hspace{0.17cm}f_{6}^{0}(x)\right),
\end{equation}
subject to 
$$
 S=\left\{x\in\mathbb{R}^{2}\hspace{0.15cm}:\hspace{0.15cm}Ax\leq b,\hspace{0.15cm} x\geq 0,\hspace{0.15cm} b\in\mathbb{R}^{9}\right\},
$$
where
$$
f_{k}^{0}(x)=c_{0k}^{t}x,\hspace{0.4cm}\text{for}\hspace{0.14cm}k\in\mathcal{K}=\left\{1,\ldots,6\right\},
$$
and
$$
\begin{array}{ccc}
A=\left(\begin{array}{rr}
1& -1\\
2& -1\\
 2& 1\\
 1& 2\\
0&1\\
 -1&1\\
 1&0\\
 2&1\\
 3&3
\end{array}\right),&\hspace{0.3cm}
b=\left(\begin{array}{r}3\\ 9\\ 19\\20 \\8\\6\\0\\4\\7
\end{array}\right),&
\begin{array}{lll} 
c_{01}=\left(\begin{array}{r} 
\frac{4}{3} \\
-1
     \end{array}\right),&c_{02}=\left(\begin{array}{r} 
\frac{4}{3}\vspace{0.15cm}\\
-2
     \end{array}\right),&c_{03}=\left(\begin{array}{r} 
\frac{5}{4}  \vspace{0.15cm}\\
-\frac{3}{4}
     \end{array}\right),\vspace{0.35cm}\\
     c_{04}=\left(\begin{array}{r} 
6\\
0
\end{array}\right),&c_{05}=\left(\begin{array}{r} 
1 \\
-2
     \end{array}\right),&c_{06}=\left(\begin{array}{r} 
1  \\
4
     \end{array}\right).\\
     
     \end{array}
\end{array}
$$
Then, generate the set \(\mathcal{A}(F^{0}) \subseteq \mathcal{L}\) of all linear mappings \(G \in \mathcal{L}\) that satisfy:
\begin{equation}\label{sanex}
   \mathcal{S}(F^{0}) = \mathcal{S}(G), \quad \text{for all} \quad G \in \mathcal{A}(F^{0}),
\end{equation}
The set \(\mathcal{V}\left(F^{0}\right) \) of vertices, ordered from \( v^{1} \) to \( v^{9} \) in a counterclockwise direction, is given by
$$
\mathcal{V}(S) = \left\{v^{i} \hspace{0.2cm}:\hspace{0.2cm} j \in V(S) = \left\{1, \ldots, 9\right\}\right\},
$$
with
\[ 
\begin{array}{cccccc}
   v^{1} = \left(4, 1\right), & v^{2} = \left(6, 3\right), & v^{3} = \left(7, 5\right), & v^{4} = \left(6, 7\right), & v^{5} = \left(4, 8\right), \vspace{0.3cm} \\
   v^{6} = \left(2, 8\right), & v^{7} = \left(0, 6\right), & v^{8} = \left(0, 4\right), & v^{9} = \left(1, 2\right). & \\
\end{array}
\]

\begin{figure}[htbp]
    \centering
     \begin{minipage}{0.43\textwidth}
        \centering        \includegraphics[width=6.5cm,height=4.5cm]{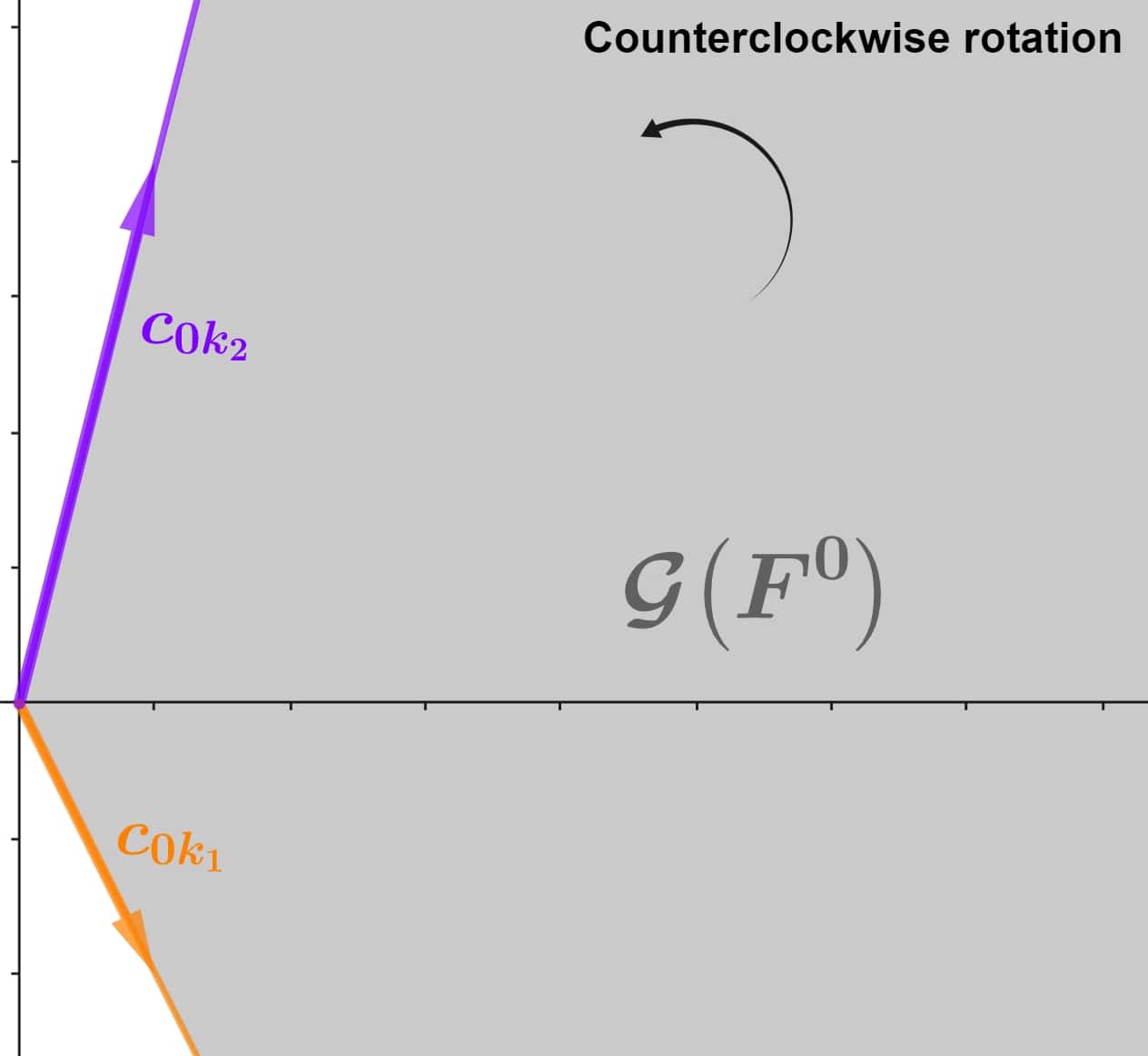}
        \caption{
The gradient cone \(\mathcal{G}(F^{0})\) and its two generators \(c_{0k_{1}}\) and \(c_{0k_{2}}\) for \(k_{1}=5\) and \(k_{2}=6\).\label{gco} }
    
    \end{minipage}\hfill
    \begin{minipage}{0.43\textwidth}
        \centering
        \includegraphics[width=7cm,height=5cm]{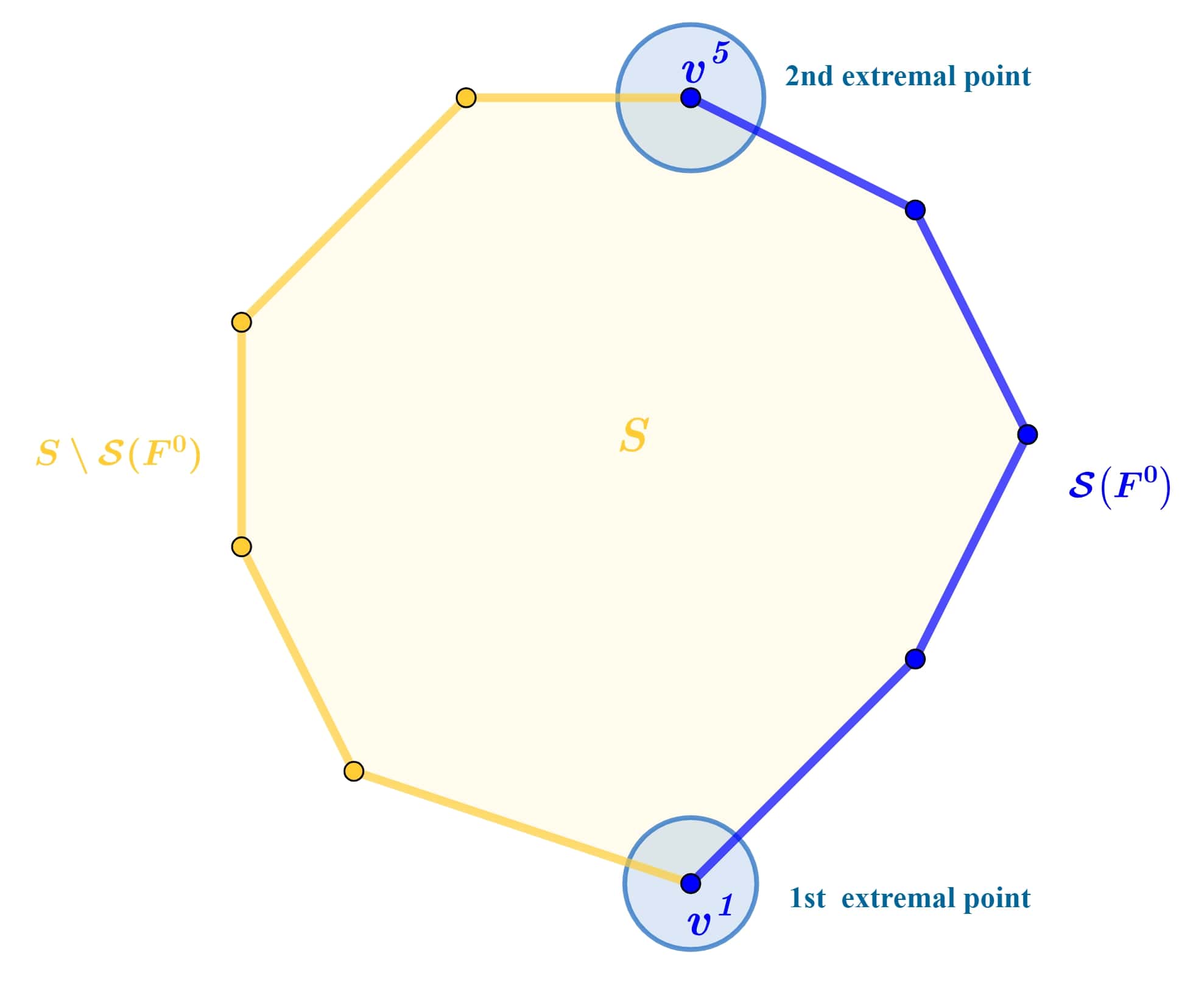}
        \caption{The two extremal efficient solutions \(v^{1}\) and  \(v^{5}\) of the set $\mathcal{S}(F^{0})$.\label{exp}}
    \end{minipage}

\end{figure}

    \begin{description}
\item[Part 1: Generate the efficient points set \( \mathcal{S}\left(F^{0}\right) \)  of problem (\ref{exemple2}).] $\;$\\
\begin{description}
        \item[Step 1:]  Express \(c_{0k}\) for all \(k \in\mathcal{K}\) in polar coordinates, as follows:
\[
\begin{array}{ll}
c_{01}=\frac{5}{3} \left(\cos(-36.870^\circ), \sin(-36.870^\circ)\right), & 
c_{02}=\frac{2\sqrt{13}}{3} \left(\cos(-56.310^\circ), \sin(-56.310^\circ)\right), \vspace{0.23cm} \\

c_{03}=\frac{\sqrt{34}}{4} \left(\cos(-30.964^\circ), \sin(-30.964^\circ)\right), &
c_{04}=6\left(\cos(0^\circ), \sin(0^\circ)\right), \vspace{0.23cm} \\

c_{05}=\sqrt{5} \left(\cos(-63.435^\circ), \sin(-63.435^\circ)\right), &
c_{06}=\sqrt{17} \left(\cos(75.964^\circ), \sin(75.964^\circ)\right)
\end{array}
\]
Denote the angles as follows:
\[
\begin{array}{llllll}
    \phi(1)=-36.87^\circ, & \phi(2)=-56.31^\circ, & \phi(3)=-30.964^\circ,&
    \phi(4)=0^\circ, &
    \phi(5)=-63.435^\circ, &
    \phi(6)=75.964^\circ.
\end{array}
\]
\item[Step 2:] The two gradients \(c_{05}\) and \(c_{06}\) generate the two extreme rays of \(\mathcal{G}(F^{0})\) (see Figure \ref{exp}). Indeed, we have:
\[
\phi(5) = \min_{k \in \mathcal{K}} \phi(k) \quad \text{and} \quad \phi(6) = \max_{k \in \mathcal{K}} \phi(k).
\]

\item[Step 3:] Using Corollary \ref{re}, we get the two extremal efficient solutions  
\(v^{1}=(4,1)\) and  \(v^{5}=(4,8)\), see Figure \ref{exp}. 
Then,
\[
\mathcal{VS}(F^{0}) = \{v^{j} \; : \; j = 1, \ldots, 5\}.
\]

\item[Step 4:] Use formula (\ref{cur}) by taking \(j^{0} = 1\) and \(j = 5\), then obtain \(\mathcal{S}(F^{0})\) (see Figure \ref{SF}) as follows:
\[
\mathcal{S}(F^{0}) := \mathcal{S}_{1}^{5} = \mathcal{H}\{v^{1}, v^{2}\} \cup \mathcal{H}\{v^{2}, v^{3}\} \cup \mathcal{H}\{v^{3}, v^{4}\} \cup \mathcal{H}\{v^{4}, v^{5}\}.
\]
\end{description} 
\begin{figure}[htbp]
    \centering
    \begin{minipage}{0.3\textwidth}
        \centering
        \includegraphics[width=6.5cm,height=7cm]{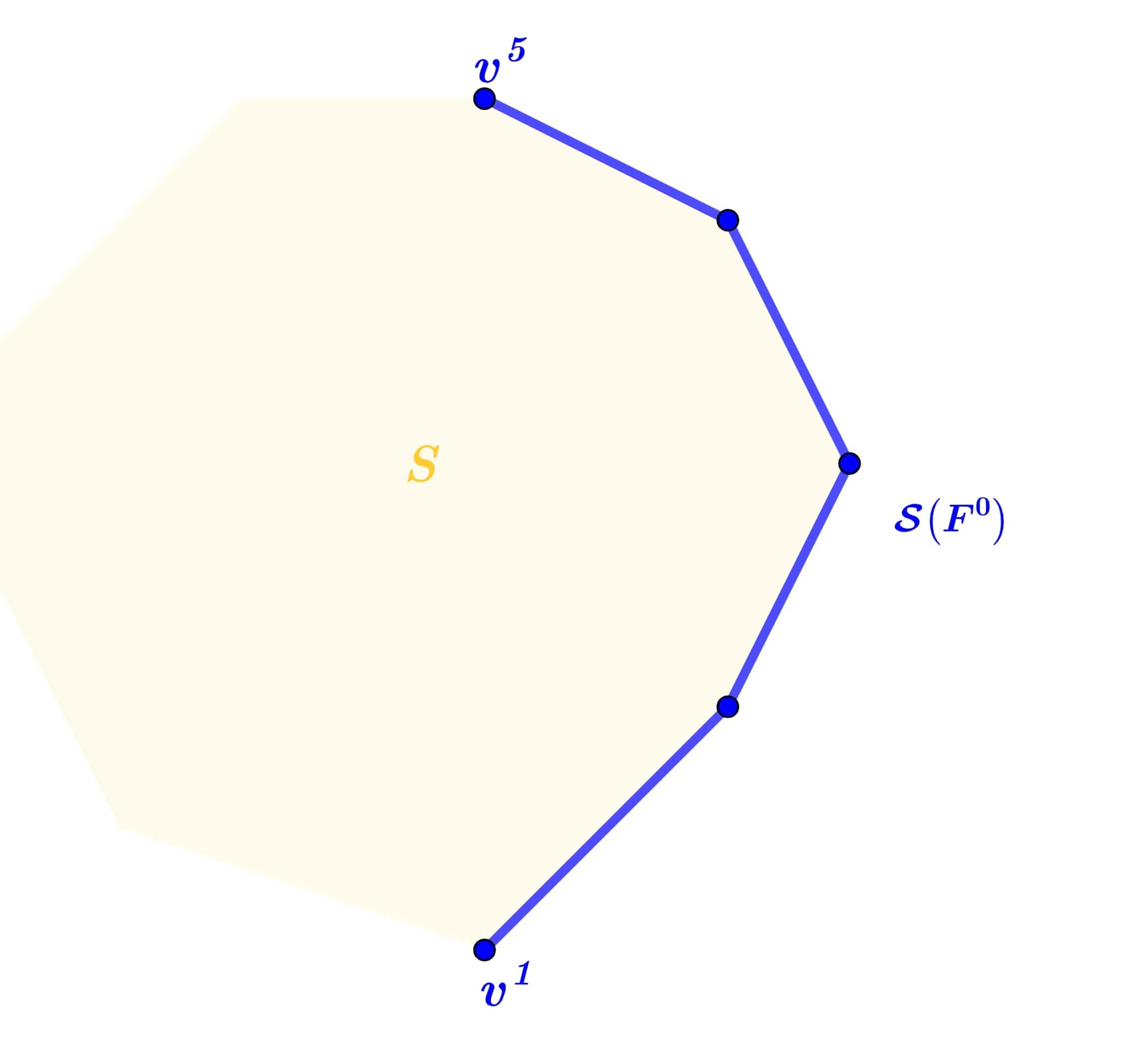}
        \caption{The efficient solutions set \( \mathcal{S}\left(F^{0}\right) \)  of problem (\ref{exemple2})\label{SF}}
    \end{minipage}\hfill
    \begin{minipage}{0.3\textwidth}
        \centering
        \includegraphics[width=3cm,height=3cm]{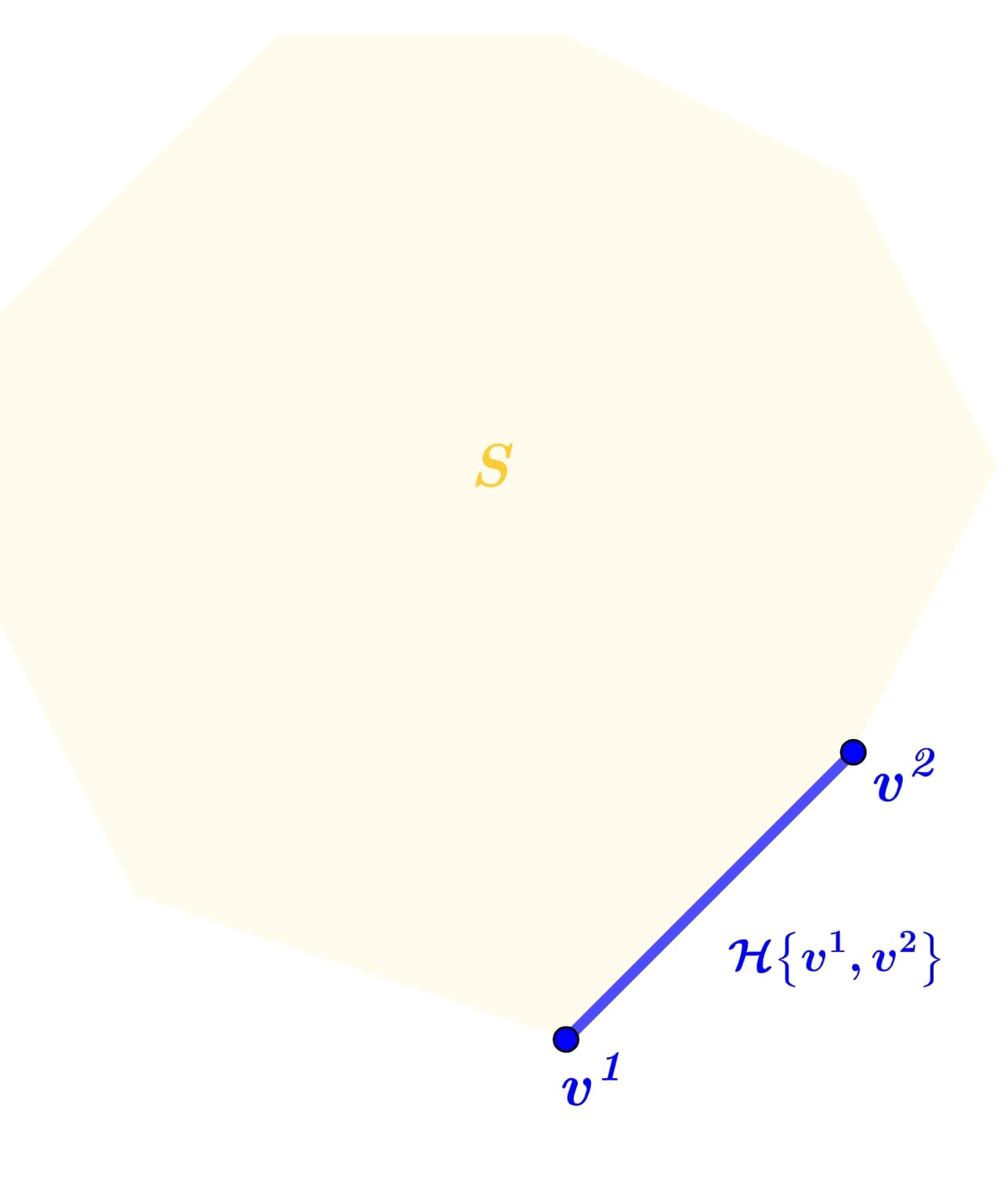}
        \vspace{0.3cm} 
        \includegraphics[width=3cm,height=3cm]{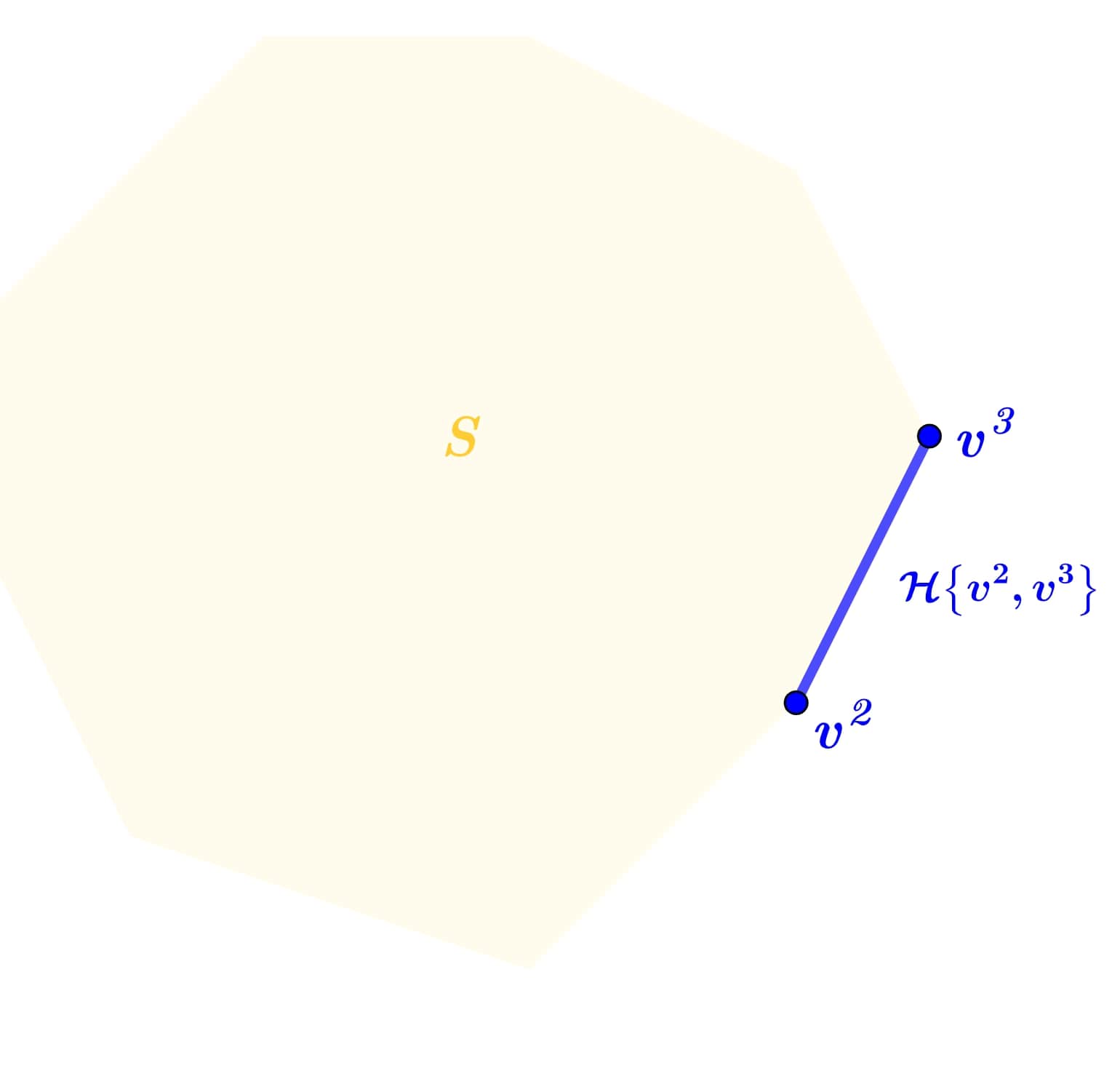}
    \end{minipage}\hfill
    \begin{minipage}{0.3\textwidth}
        \centering
        \includegraphics[width=3cm,height=3cm]{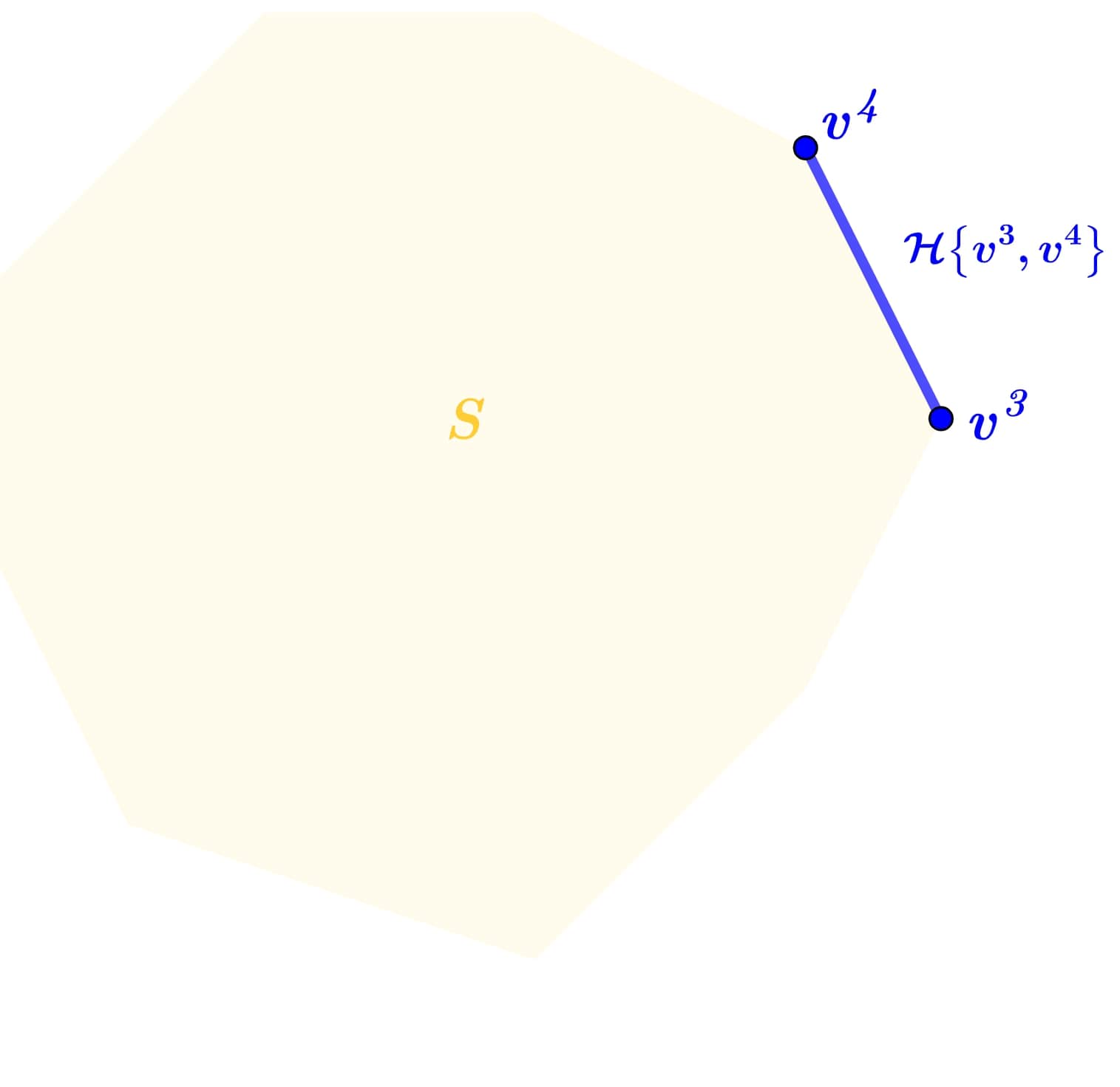}
        \vspace{0.3cm} 
        \includegraphics[width=3cm,height=3cm]{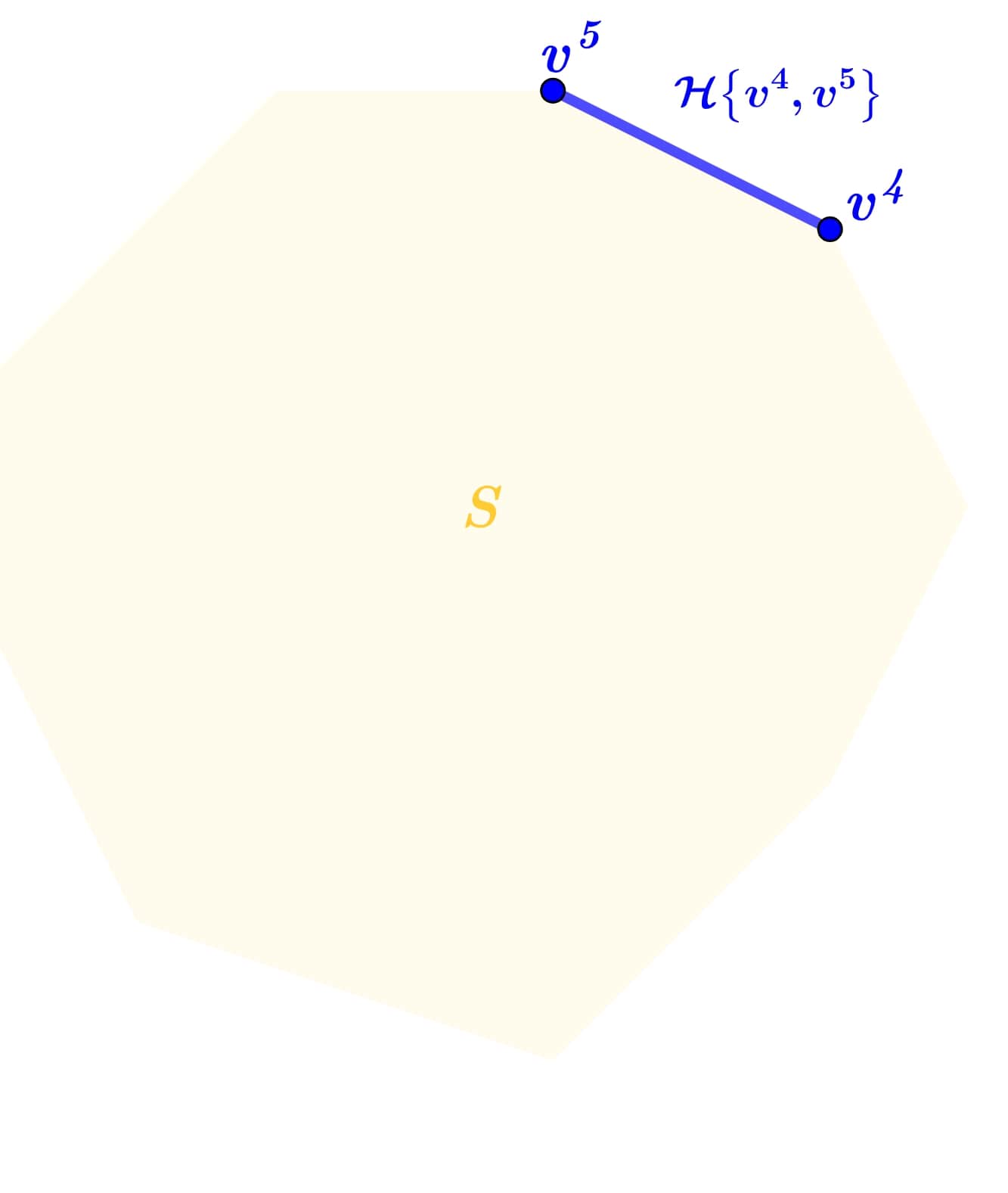}
    \end{minipage}
\end{figure}
\item[Part 2: Generate the solution set \(\mathcal{A}(F^{0})\) for the sensitivity analysis problem (\ref{sanex}).] $\;$\\
\begin{description}
\item[Step 5:] By formulas (\ref{pok}), we get \(\theta_{1}(5) = -18.43^\circ\) and \(\theta_{2}(6) = 180^\circ\). Indeed, we have:
\[
\begin{array}{lcccl}
   v^{12} &=& v^{1} - v^{9} &=& \sqrt{10}\left(\cos(-18.435^\circ), \sin(-18.435^\circ)\right), \vspace{0.23cm}\\
   v^{56} &=& v^{6} - v^{5} &=& 2\left(\cos(180^\circ), \sin(180^\circ)\right).
\end{array} 
\]
    \item[Step 6:]  By Corollary \ref{ffd}, we obtain  \(\mathcal{A}(F^{0})\) as follows:

\[
\mathcal{A}(F^{0}) = \bigcup_{K \geq 2} \left( \bigcup_{(\theta_{1}, \ldots, \theta_{K-2}) \in \hat{I}(k_{1},k_{2})^{K-2}} \left\{ \left( \langle c_{05}^{t}, \cdot \rangle, \langle c_{06}^{t}, \cdot \rangle, \langle R_{\theta_{1}}(c_{05})^t, \cdot \rangle, \ldots, \langle R_{\theta_{K-2}}(c_{05})^t, \cdot \rangle \right) \right\} \right),
\]
where
\[
\hat{I}(k_{1},k_{2})^{K-2} = \left]-108.43^\circ, 90^\circ\right[.
\]

    \end{description}
    \end{description}

\section{Discussion}\label{section 6}

The proposed method has been introduced as a sensitivity analysis approach for the set of efficient solutions of a MOLP problem. Additionally, it can be regarded as a classification method for MOLP problems (see Remark \ref{CLA}).

By conceptualizing objects as geometric entities, objective functions as graphs, and feasible regions as convex sets, we have identified two significant components. First, the weighted function of the objective function \( F_{\lambda}^{0} \) can be interpreted as a rotation of one of its components \( f_{k_{1}}^{0} \). Specifically, assigning a weight vector \(\lambda\) to the objective function \( F^{0} \) results in rotating the graph of its component \( f_{k_{1}}^{0} \) by an angle \(\theta\) (Theorem \ref{lem}). The second component is Theorem \ref{conn}, which asserts that the efficient solutions are connected; this means one cannot transition from one efficient point to another without passing through a sequence of adjacent efficient points.

Utilizing the fact that \( S \) is within \(\mathbb{R}^{2}\), and adopting the convention that rotation is counterclockwise (see Remark \ref{sr}), we initially demonstrate that the extreme points associated with the two generators \( c_{0k_{1}} \) and \( c_{0k_{2}} \) of the two extreme rays of the gradient cone \(\mathcal{G}(F^{0})\) are extremal points (see Corollary \ref{re} and Figure \ref{exp}). Since transitioning from the first extreme ray can be done through both clockwise and counterclockwise paths, which generate two different sets of efficient solutions, Remark \ref{sr} is essential for determining the direction of traversal. This foundation allows us to prove in Theorem \ref{fff} the equivalence between MOLP problems and TOLP problems. In \(\mathbb{R}^{2}\), all MOLP problems can be reduced to TOLP problems, meaning that for each MOLP problem, there exists an equivalent TOLP problem with the same set of efficient solutions, and vice versa.

Remark \ref{CLA} provides a classification of MOLP problems by associating each element of \(\mathcal{NS}\) with a subset of \(\mathcal{L}\) in a bijective manner. It is important to note that the idea of constructing \(\mathcal{NS}\) stems from the fact that the set of efficient solutions forms a continuous curve. Thus, the idea was to consider all finite successive sequences of elements from \(\mathcal{V}(S)\) and to consider the polygonal curve connecting the elements of each sequence.

Algorithm \ref{ama} serves as a comprehensive guide for applying the approach developed in this work, demonstrating the effectiveness of the method. The first part consists of the initial four steps, which generate the set of efficient solutions \(\mathcal{S}(F^{0})\). We begin by expressing the gradients of the objective functions in polar coordinates and then identify the generators \( c_{0k_{1}} \) and \( c_{0k_{2}} \) of the extreme rays of the gradient cone by examining the maximum and minimum angles. This simplifies the identification of the extreme rays. Next, we maximize \( f^{0}_{k_{1}} \) and \( f^{0}_{k_{2}} \) to obtain the extremal points \( v^{j_{k_{1}}} \) and \( v^{j_{k_{2}}} \). In step 4, we take the union of the convex combinations of successive extreme points from \( v^{j_{k_{1}}} \) to \( v^{j_{k_{2}}} \) to form the polygonal curve \(\mathcal{S}(F^{0})\). This demonstrates the simplicity of the calculations required to generate efficient solutions compared to traditional procedures.\\
The second part involves two straightforward steps that do not require extensive calculations to generate the set \(\mathcal{A}(F^{0})\). In the fifth step, we convert the two vectors into polar coordinates as described in formulas (\ref{pok}), identifying the angles \(\theta_{1}(k_{1})\) and \(\theta_{2}(k_{2})\). Finally, in the sixth step, we apply formula (\ref{sets}) to construct the desired set.

\section{Conclusion}\label{section 7}

In this work, we have introduced a novel geometric approach for the sensitivity analysis of MOLP problems. This method builds upon the geometric sensitivity analysis approach provided by Kaci and Radjef \cite{kaci2022}. We began by visualizing the weighting of the objective function as a rotation of the graph of one of its components. We then demonstrated that this rotation occurs between the two extreme rays of the gradient cone, generating all efficient solutions. This significantly simplifies calculations and allows us to prove the equivalence between MOLP and TOLP problems, thereby facilitating their classification.

The classification of MOLP problems based on their sets of efficient solutions has not been addressed in the existing literature. This work pioneers such a classification approach and introduces a novel method for performing sensitivity analysis on MOLP problems. This innovative perspective explains the absence of numerical comparisons in this paper.

Graphical illustrations accompany the technical details to aid in understanding and visualizing the results. In particular, Figure \ref{GH} is designed to be consulted progressively while reading Section \ref{section 4}. Algorithm \ref{ama} connects all the results and serves as a practical guide for applying them. A numerical example and a discussion of the results are also provided.

 \bibliographystyle{acm}
\bibliography{references}

\end{document}

%% file: parameters.tex
\usepackage{setspace}
\usepackage{lmodern}

\usepackage{textgreek}

\usepackage{fixltx2e}

\reversemarginpar

\usepackage{amsmath,amssymb,amsfonts,amsthm}
\usepackage{mathrsfs}

\usepackage[locale = FR]{siunitx} 
\DeclareSIUnit\vitesse{\meter\per\second}
\usepackage{eurosym}
\DeclareSIUnit{\octet}{o}

\usepackage{graphicx,array,tikz,multirow}
\usepackage{caption,subcaption} 
\usepackage{svg,float}
\usepackage{booktabs,paralist}
\newcolumntype{x}[1]{>{\centering\arraybackslash\hspace{0pt}}p{#1}}
\usepackage[section]{placeins}
\usepackage{hanging}

\usepackage{fancyhdr,emptypage} 

\fancypagestyle{plain}{ 
    \fancyhead{}\fancyfoot[C]{\thepage}

}

\usepackage[Conny]{fncychap}

\usepackage[francais,nohints,tight]{minitoc}		

\usepackage[nottoc]{tocbibind} 
\usepackage[titles]{tocloft}

\usepackage{textcomp} 

\usepackage{titling}
\usepackage{lipsum} 
\usepackage{csquotes} 

\usepackage{xspace}
\usepackage{afterpage}

\usepackage[textsize=footnotesize]{todonotes}

\usepackage{bookmark}
\usepackage{acronym}
\usepackage[nameinlink,french]{cleveref} 
\Crefname{figure}{Fig.}{Figs.} 
\crefname{figure}{fig.}{figs.}
\Crefname{equation}{Eq.}{Eqs.}
\crefname{equation}{eq.}{eqs.}
\Crefname{table}{Table.}{Tables.}
\crefname{table}{table.}{tables.}

\definecolor{color_ref}{rgb}{1.0, 0.13, 0.32} 
\definecolor{color_link}{rgb}{0.0, 0.0, 1.0}
\definecolor{curcolor}{rgb}{0.0, 0.0, 1.0} 
\definecolor{brightpink}{rgb}{1.0, 0.0, 0.5} 
\definecolor{navyblue}{rgb}{0.0, 0.0, 1.0}
\hypersetup{
	colorlinks=true, 
	pdfstartview=FitV, 
	urlcolor=brightpink, 
	linkcolor= navyblue, 
	citecolor=color_ref 
}
